# SUR L'INTERSECTION DES COURANTS LAMINAIRES

ROMAIN DUJARDIN

ABSTRACT. We try to find a geometric interpretation of the wedge product of positive closed laminar currents in $\mathbb{C}^2$. We say such a wedge product is *geometric* if it is given by intersecting the disks filling up the currents.

Uniformly laminar currents do always intersect geometrically in this sense. We also introduce a class of "strongly approximable" laminar currents, natural from the dynamical point of view, and prove that such currents intersect geometrically provided they have continuous potentials.

## 1. INTRODUCTION

Une motivation pour l'étude du produit extérieur des courants positifs fermés au cours des années 70 était de donner un nouveau point de vue sur l'intersection des ensembles analytiques. Il apparaît que dans le cas de courants d'intégration sur des diviseurs, l'intersection pluripotentialiste et l'intersection algébro-géométrique coïncident. Ce sont des raisons analogues qui ont poussé E. Bedford, M. Lyubich et J. Smillie (voir [BLS] pp. 79-80) à introduire les courants laminaires, et à étudier leur intersection potentialiste sous un angle géométrique, s'appuyant en cela sur une construction de D. Ruelle et D. Sullivan [RS]. Il y est effectivement montré que dans le cas des automorphismes polynomiaux de $\mathbb{C}^2$, les courants dynamiques stable et instable ont une structure laminaire et que leur produit extérieur admet une interprétation géométrique ; la preuve repose de façon essentielle sur la théorie de Pesin des systèmes dynamiques non uniformément hyperboliques. Ce résultat d'intersection géométrique est utilisé pour montrer que les points périodiques de type selle sont situés dans le support de la mesure d'entropie maximale et engendrent des intersections homoclines.

Nous allons étudier ces questions sous un angle plus systématique, et retrouver ces résultats dans un cadre indépendant de toute dynamique. Ceci fournit en contrepartie une preuve des résultats de [BLS] mentionnés ci-dessus, qui n'utilise pas la théorie de Pesin.

Rappelons qu'un courant positif fermé de bidegré (1,1) dans un ouvert $\Omega \subset \mathbb{C}^2$ est dit *uniformément laminaire* s'il est localement un courant d'intégration sur une lamination transversalement mesurée et *laminaire* s'il est localement défini comme intégrale directe d'une famille de disques compatibles (voir la section 2 pour plus de détails). On dira que le produit extérieur de deux courants laminaires est *géométrique* s'il est décrit par les intersections des disques constituant ces courants, avec la convention que l'intersection d'un disque avec lui même est nulle. En particulier l'auto-intersection d'un courant laminaire est géométrique si et seulement si elle est nulle.

On ne considèrera que des courants diffus, la théorie de l'intersection des courants d'intégration étant bien comprise, et tout courant positif fermé se décomposant d'après le théorème de semi-continuité de Siu en une somme de courants d'intégration et d'un courant ne chargeant pas les sous ensembles analytiques (voir Demailly [De1]).

Nous montrons à la section 3 que le produit extérieur des courants uniformément laminaires est toujours géométrique –pourvu qu'il soit défini. Des exemples simples montrent qu'il n'est pas raisonnable d'espérer un résultat aussi général pour les courants laminaires. Les courants dynamiques apparaissant en dynamique des automorphismes polynomiaux de $\mathbb{C}^2$ (et plus généralement, des applications birationnelles) ont toutefois une propriété supplémentaire : ils sont limites de suites de diviseurs rationnels $\frac{1}{d_n}[C_n]$ de topologie contrôlée. Nous dirons de ces courants qu'ils sont *fortement approximables* (voir la section 4 pour une définition formelle). Notre résultat principal est le suivant :

**Théorème.** *Soient $T_1$ et $T_2$ des courants laminaires fortement approximables dans $\Omega$, de potentiel continu. Alors le produit extérieur $T_1 \wedge T_2$ est géométrique dans $\Omega$.*

Un corollaire de ce théorème est que tous les courants laminaires ne sont pas fortement approximables : il est aisé de trouver des exemples de courants laminaires de potentiel continu satisfaisant $T \wedge T > 0$. Ceci suggère que la classe des courants laminaires fortement approximables exhibe des propriétés intéressantes. Nous prévoyons de les étudier plus avant dans un article à venir.

Le plan de l'article est le suivant. Nous rappelons dans le premier paragraphe quelques faits concernant les courants laminaires et leur intersection. À la section 3 nous étudions l'intersection des courants uniformément laminaires, et à la section 4 celle des courants fortement approximables.

Cet article est issu de la thèse de doctorat [Du3] de l'auteur.

## 2. Définitions et quelques exemples

Nous commençons par rappeler brièvement les diverses notions de laminarité [BLS, Ca, Du1, Du3] pour les courants positifs fermés. On fixe un ouvert $\Omega \subset \mathbb{C}^2$. $T$ étant un courant positif, on notera Supp($T$) son support, $\mathbf{M}(T)$ sa masse, et $\|T\|$ sa mesure trace. La notation $\mathbb{D}$ désigne le disque unité dans $\mathbb{C}$.

**Définition 2.1.**
*i. Le courant $T$ est uniformément laminaire si pour tout $x \in$ Supp($T$), il existe des ouverts $V \supset U \ni x$, tels que $V$ soit biholomorphe au bidisque $\mathbb{D}^2$ et dans la carte correspondante de $V$, $T|_U$ soit l'intégrale directe d'une famille de courants d'intégration sur des graphes disjoints dans $\mathbb{D}^2$, plus précisément :*

*il existe une mesure finie $\mu$ sur $\{0\} \times \mathbb{D}$ et une famille de graphes disjoints $(\Gamma_a)$ dans $\mathbb{D}^2$ tels que $(0, a) \in \Gamma_a$, et*

$$(1) \qquad T|_U = \int_{\{0\} \times \mathbb{D}} [\Gamma_a \cap U] d\mu(a).$$

*ii. Le courant $T$ est laminaire dans $\Omega$ s'il existe une suite d'ouverts $\Omega^i \subset \Omega$ tels que $\|T\|(\partial\Omega^i) = 0$ et une suite croissante de courants $(T^i)_{i\geq 0}$, $T^i$ uniformément laminaire dans $\Omega^i$, telles que*

$$\lim_{i\to\infty} T^i = T.$$

Un courant uniformément laminaire est un cycle feuilleté sur une lamination par surfaces de Riemann, plongée et munie d'une mesure transverse invariante. On montre facilement (en utilisant [BLS]) que $T$ est laminaire au sens de la définition 2.1 si et seulement s'il est faiblement laminaire au sens de [BLS], c'est à dire que $T$ s'écrit localement comme intégrale d'une famille mesurée de disques holomorphes *compatibles*, i.e. on a la représentation

$$(2) \qquad\qquad T = \int_{\mathcal{A}} [D_a] d\mu(a),$$

où deux disques s'intersectant ont un ouvert en commun. Dans les cas que nous considèrerons par la suite, $\Omega^i = \Omega\backslash\partial\mathcal{Q}_i$ où $\mathcal{Q}_i$ est une subdivision de $\Omega$ à un ensemble de $\|T\|$ mesure nulle près. Il est important de noter que la représentation d'un courant laminaire –comme limite croissante ou intégrale d'une famille de disques– n'est jamais unique puisqu'il est toujours possible de la modifier sur un ensemble de $\|T\|$ mesure nulle. On dira qu'un courant est *diffus* s'il ne charge aucun sous ensemble analytique.

Soient $T_1$ et $T_2$ deux courants positifs fermés dans $\Omega$. Supposons que le courant $T_1$ admette un potentiel $u_1$ dans $\Omega$, i.e. $T_1 = dd^c u_1$ où $dd^c = \frac{i}{\pi}\partial\bar\partial$. On dira que le produit extérieur $T_1 \wedge T_2$ est *admissible* si $u_1 \in L^1_{loc}(\|T_2\|)$. Ceci ne dépend clairement pas du choix de $u_1$. On définit alors le produit extérieur $T_1 \wedge T_2$ par la formule usuelle

$$T_1 \wedge T_2 = dd^c(u_1 T_2)$$

qui est une mesure positive de masse localement finie. N. Sibony a montré (voir [Du3]) que le produit extérieur ainsi défini est symétrique –autrement dit $T_1 \wedge T_2$ admissible implique $T_2 \wedge T_1$ admissible et $T_1 \wedge T_2 = T_2 \wedge T_1$– et satisfait un théorème de convergence décroissante. On voit également que si $T_1' \leq T_1$ et $T_2' \leq T_2$ sont des courants positifs fermés alors

$$(3) \qquad\qquad T_1 \wedge T_2 \text{ admissible } \Rightarrow T_1' \wedge T_2' \text{ admissible et } T_1' \wedge T_2' \leq T_1 \wedge T_2.$$

En effet on peut localement choisir un potentiel $u_i$ de $T_i$ (resp. $u_i'$ de $T_i'$), $i = 1, 2$, tels que $u_i' \geq u_i$.

Cette notion nous semble bien adaptée au produit extérieur des courants (uniformément) laminaires, notamment parce qu'elle est stable par restriction du domaine. On aurait cependant pu également adopter la condition d'admissibilité de [FS], qui satisfait également la relation (3).

Étant donnés deux disques holomorphes $\Delta_1$ et $\Delta_2$, définissons la mesure $[\Delta_1 \cap \Delta_2]$ comme la somme des masses de Dirac aux points d'intersection de $\Delta_1$ et $\Delta_2$, comptés sans multiplicité (celle ci n'a pas d'influence sur l'intersection des courants diffus, et on aurait tout autant pu la compter) si ceux-ci sont isolés, $[\Delta_1 \cap \Delta_2] = 0$ si les deux disques ont un ouvert en commun.

Si $T_1$ et $T_2$ sont des courants *uniformément laminaires* diffus, s'écrivant dans l'ouvert $U$ sous la forme

$$T_i = \int_{\mathcal{A}_i} [\Delta_{i,a}] d\mu_i(a), \; i = 1, 2,$$

(où les $\Delta_{i,a}$ sont des sous variétés de $U$), et admettant un produit extérieur, on définit la *mesure d'intersection géométrique* dans $U$ comme

$$T_1 \dot{\wedge} T_2 = \int_{\mathcal{A}_1 \times \mathcal{A}_2} [\Delta_{1,a} \cap \Delta_{2,b}] d\mu_1(a) d\mu_2(b),$$

(en particulier $T_1 \dot{\wedge} T_1 = 0$) et on dira que l'intersection de $T_1$ et $T_2$ est *géométrique* dans $U$ si $T_1 \wedge T_2 = T_1 \dot{\wedge} T_2$. Nous verrons à la section suivante que c'est en fait toujours le cas.

Soient maintenant des courants *laminaires* fermés diffus $T_1$ et $T_2$ admettant un produit extérieur. $T_1$ et $T_2$ admettent des représentations comme limites croissantes $T_k = \lim T_k^i$, avec $T_k^i$ uniformément laminaire dans $\Omega_k^i$ ; nous dirons que ces représentations sont *compatibles* si pour tout $i$, $\partial \Omega_k^i$ est de $\|T_j\|$ mesure nulle, $j, k \in \{1, 2\}$. On dira que le produit extérieur de $T_1$ et $T_2$ est *géométrique* s'il existe des représentations de $T_1$ et $T_2$ comme limites croissantes comme à la définition 2.1 telles que la suite des mesures d'intersection géométrique $T_1^i \dot{\wedge} T_2^i$ (définie dans $\Omega_1^i \cap \Omega_2^i$) croisse vers $T_1 \wedge T_2$. En particulier l'auto-intersection d'un courant laminaire est géométrique si et seulement si elle est nulle.

En termes de représentation comme intégrales de disques comme en (2), l'intersection est géométrique si et seulement si on a l'égalité

$$T_1 \wedge T_2 = \int_{\mathcal{A}_1 \times \mathcal{A}_2} [D_a^1 \cap D_b^2] d\mu_1(a) d\mu_2(b).$$

Nous appellerons *mesure d'intersection géométrique*, le membre de droite de cette égalité, ou de façon équivalente la limite croissante des mesures $T_1^i \dot{\wedge} T_2^i$.

On voit que la mesure d'intersection géométrique est très sensible à la représentation des courants. En particulier si $T_1 \wedge T_2$ est singulière par rapport à $\|T_1\|$ et $\|T_2\|$, il est toujours possible de modifier les représentations de manière que la mesure d'intersection géométrique associée soit nulle. Il est donc essentiel de pouvoir disposer de représentations adaptées des courants laminaires. Ce sera une motivation pour l'introduction des courants fortement approximables (définition 4.1 ci-après).

Les exemples suivants illustrent la difficulté de l'interprétation géométrique de l'intersection des courants laminaires généraux.

**Exemple 2.2.** Cet exemple est une manifestation du phénomène évoqué ci dessus. Soit $X$ un compact non polaire et de mesure nulle dans $\mathbb{D}$, et $\mu$ une mesure à potentiel continu portée par $X$. Soient $T^h$ et $T^v$ deux courants respectivement de supports horizontal et vertical dans $\mathbb{D}^2$, définis par

$$T^h = \int_X [\mathbb{D} \times \{w\}] d\mu(w) \text{ et } T^v = \int_X [\{z\} \times \mathbb{D}] d\mu(z).$$

Bien entendu le produit $T^h \wedge T^v = \mu \otimes \mu$ est géométrique. Cependant, si on subdivise chacun des $\{z\} \times (\mathbb{D} \backslash X)$ et $(\mathbb{D} \backslash X) \times \{w\}$ en une réunion de disques disjoints, on obtient des représentations de $T^h$ et $T^v$ comme courants laminaires et la mesure d'intersection géométrique associée est nulle, puisque les disques ainsi construits ne se coupent pas.

De même le courant $T = T^h + T^v$ est laminaire, de potentiel continu et $T \wedge T > 0$. Cette auto-intersection n'est donc pas géométrique, et illustre le fait que des disques subordonnés au courant s'intersectent ; ceci n'est pas le cas dans l'exemple suivant.

**Exemple 2.3.** J.P. Demailly [De2] a introduit le courant positif fermé dans $\mathbb{C}^2$

$$T = dd^c \max(\log^+ |z|, \log^+ |w|)$$

qui se décompose de la façon suivante :

$$T = \int_{S^1} [\{e^{i\theta}\} \times \mathbb{D}] d\lambda(\theta) + \int_{S^1} [\mathbb{D} \times \{e^{i\theta}\}] d\lambda(\theta) + \int_{S^1} [V_\theta] d\lambda(\theta),$$

où $\lambda$ est la mesure de Lebesgue sur le cercle unité $S^1$, et $V_\theta = \{(z,w) \in \mathbb{C}^2, \ z = e^{i\theta}w, |z| > 1\}$. Le courant $T$ est donc un courant laminaire, fermé, et de potentiel continu. Le courant $T$ se prolonge par ailleurs comme courant positif fermé de potentiel continu dans le plan projectif $\mathbb{P}^2$. On montre aisément que $T \wedge T > 0$ est la mesure de Lebesgue sur le tore unité. Cette auto-intersection n'est donc pas géométrique.

**Exemple 2.4.** L'auto-intersection (au sens de la théorie du pluripotentiel) d'un courant laminaire ne provient pas nécessairement de l'intersection géométrique des disques ou de leurs prolongements, comme le montre l'exemple suivant. Soit un ensemble de Cantor de dimension $< 1$ dans le plan et $G_K$ sa fonction de Green. On suppose que $G_K$ est continue, et on a $G_K = 0$ sur $K$, $G_K > 0$ hors de $K$. Soit $G(z,w) = \max(G_K(z), G_K(w))$, on a

$$\text{Supp}(dd^c G) \subset \{G_K(z) = G_K(w) > 0\} \cup (K \times K).$$

Comme $\dim(K \times K) < 2$, $K \times K$ n'est pas chargé par $dd^c G$ ; de plus l'hypersurface

$$\{G_K(z) = G_K(w) > 0\} = \{G_K(z) - G_K(w) = 0\} \setminus (K \times K)$$

est Levi-plate, donc feuilletée par des courbes analytiques. On en déduit que le courant $dd^c G$ est laminaire, uniformément dans le complémentaire de $K \times K$ par le théorème de support de Demailly [De1]. De plus $(dd^c G)^2$ est une mesure de probabilité portée par $K \times K$ car $G$ est la fonction de Green pluricomplexe de $K \times K$.

Il est possible de choisir $K$ de telle sorte que les courbes feuilletant $\{G_K(z) = G_K(w) > 0\}$ ne se prolongent pas à travers $K$ (à l'exception de la droite $(z = w)$). Soit un point $(z_0, w_0)$ de $K \times K$ et un disque holomorphe local $\Delta \ni (z_0, w_0)$, que l'on pourra supposer être un graphe $\Delta = \{(\zeta, \varphi(\zeta)), \zeta \in U\}$ au dessus de la première coordonnée, tel que $\Delta \setminus (K \times K) \subset \{G_K(z) = G_K(w) > 0\}$.

On a alors $\varphi(z_0) = w_0$ et sur $U \setminus K$, $G_K(z) = G_K(\varphi(z))$, relation qui par continuité s'étend à $U$. Quitte à se restreindre à un sous ouvert de $U$ on peut supposer que $\varphi$ n'a pas de points critiques. Il suffit donc de trouver un tel compact $K$ n'ayant pas d'automorphismes conformes locaux, auquel cas cette relation implique $\varphi = id$. On modifie pour cela la construction de l'ensemble de Cantor triadique usuel de façon que la dimension de Hausdorff locale de $K$ soit une fonction strictement croissante sur $K \subset [0,1]$. Par Tsuji [T] §III.16, on peut faire en sorte que la fonction de Green $G_K$ soit continue. $\qquad\square$

## 3. Courants uniformément laminaires

Nous montrons dans cette section que le produit extérieur des courants uniformément laminaires est géométrique dès qu'il est bien défini. Nous donnons par ailleurs un critère simple assurant l'admissibilité du produit extérieur dans ce cadre. Le cas du produit de deux courants uniformément laminaires dont l'un est à potentiel continu est traité dans le lemme 8.3 de [BLS] : leur preuve est facilement étendue au cas général puisqu'elle n'utilise qu'un argument de convergence monotone.

**Théorème 3.1.** *Soient $T_1$ et $T_2$ des courants uniformément laminaires diffus, tels que le produit extérieur $T_1 \wedge T_2$ soit admissible. Alors l'intersection $T_1 \dot\wedge T_2$ est géométrique.*

En particulier si $T_1 \wedge T_1$ est un produit extérieur admissible, alors $T_1 \wedge T_1 = 0$.

**Corollaire 3.2.** *Soient $T_1$ et $T_2$ des courants laminaires, munis de représentations (compatibles) comme limites croissantes, comme à la définition 2.1, $T_k = \lim_{i \to \infty} T_k^i$, $k = 1, 2$, et tels que le produit extérieur $T_1 \wedge T_2$ soit admissible. Alors la mesure d'intersection géométrique $\nu$ de $T_1$ et $T_2$ est dominée par $T_1 \wedge T_2$. En particulier $\nu$ est de masse localement finie.*

**Preuve du corollaire :** rappelons que les représentations comme limite croissantes sont compatibles si les bords des ouverts $\Omega^i$ ne sont pas chargés par les courants. Soit

$$\nu^i = T_1^i \dot\wedge T_2^i|_{\Omega_1^i \cap \Omega_2^i}$$

la mesure d'intersection géométrique des courants uniformément laminaires $T_1^i$ et $T_2^i$. Dans $\Omega_k^i$ on a $T_k = T_k^i + R_k^i$, $k = 1, 2$, et donc dans $\Omega_1^i \cap \Omega_2^i$ on a

$$T_1 \wedge T_2 = T_1^i \wedge T_2^i + T_1^i \wedge R_2^i + R_1^i \wedge T_2^i + R_1^i \wedge R_2^i$$

où tous les produits extérieurs sont admissibles d'après la relation (3). D'après le théorème on a l'égalité $\nu^i = T_1^i \dot\wedge T_2^i$, et donc $\nu^i \le T_1 \wedge T_2$. La mesure d'intersection géométrique étant la limite croissante des $\nu^i$ on conclut. $\qquad\square$

**Preuve du théorème :** soit $x \in \mathrm{Supp}(T_1) \cap \mathrm{Supp}(T_2)$, il suffit de montrer le résultat au voisinage de $x$. Soient respectivement $\mathcal{L}_1$ et $\mathcal{L}_2$ les laminations sous-jacentes à $T_1$ et $T_2$. Fixons un ouvert $V \ni x$ biholomorphe à un bidisque, tel que la feuille $L_1(x)$ passant par $x$ soit un graphe au dessus de la première coordonnée dans $V$. Alors pour un voisinage $V'$ de $x$, et pour tout $y \in V'$, $L_1(y')$ est un graphe dans $V$. Ainsi on peut écrire $T_1 = T_1' + T_1''$, dans $V$ où $T_1'$ est le courant uniformément laminaire formé des feuilles de $T_1$ rencontrant $V'$. $T_1'$ est fermé dans $V$ et $T_1|_{V'} = T_1'|_{V'}$. Le produit extérieur $T_1' \wedge T_2$ est admissible et il suffit donc de montrer qu'il est géométrique. Sans perte de généralité on écrit dorénavant $T_1$ pour $T_1'$.

Dans des coordonnées adaptées dans $V$ on a un potentiel pour $T_1$ de la forme

$$u_1 = \int \log|y - \psi_a(x)| \, d\mu_1(a),$$

et puisque $u_1 \in L^1(\|T_2\|)$, pour presque tout $a$, $\log|y - \psi_a(x)| \in L^1(\|T_2\|)$, c'est à dire que pour le graphe $\Gamma_a$ correspondant, le produit $[\Gamma_a] \wedge T_2$ est admissible. Montrons qu'il est géométrique.

Fixons un nouveau système de coordonnées $(z, w)$ telles que le graphe $\Gamma_a$ (maintenant fixé) ait pour nouvelle équation $\{w = 0\}$. Alors pour $\lambda$ générique, quitte à réduire encore une fois $V$ et $V'$, les feuilles de $\mathcal{L}_2$ intersectant $V'$ sont des graphes au dessus de $\{w = 0\}$ dans le système de coordonnées $(z + \lambda w, w) = (\zeta, w)$.

On écrit $T_2 = \int_{\mathcal{B}} [D_b] d\mu_2(b)$ dans $V$. Notons tout d'abord que si pour un certain $D_b$ il y a un point d'intersection non transverse dans $\Gamma_a \cap D_b$, il sera perturbé en $k$ points d'intersection transverses en changeant légèrement $b$ ([BLS], lemme 6.4). Ceci signifie que quitte à retirer une quantité au plus dénombrable de paramètres $b \in \mathcal{B}$ (ce qui ne change pas $T_2$ puisque la mesure $\mu_2$ est diffuse) on peut supposer toutes les intersections transverses.

On suit maintenant [BLS] Lemme 8.3. On écrit l'équation de $D_b$ sous la forme $w = \varphi_b(\zeta)$ et on a

$$\log|w - \varphi_b(\zeta)| = \sum_{j=1}^{N_b} \log|\zeta - p_j| + h_b(\zeta)$$

où $h_b$ est harmonique. Soient

$$\mathcal{B}_R = \{b \in \mathcal{B}, \ \|h_b\|_{L^\infty} \leq R, \ N_b \leq R\},$$

et $u_{2,R} = \int_{\mathcal{B}_R} \log|w - \varphi_b(\zeta)| \, d\mu_2(b)$. On peut supposer que tous les logarithmes sont négatifs. Quand $R \to \infty$, $\mathcal{B}_R$ croît vers $\mathcal{B}$ et $u_{2,R}$ décroît vers $u_2$. D'après le théorème de convergence monotone (standard) $u_{2,R}[\Gamma_a] \rightharpoonup u_2[\Gamma_a]$ et donc $dd^c u_{2,R} \wedge [\Gamma_a] \rightharpoonup dd^c u_2 \wedge [\Gamma_a]$. Par ailleurs le produit extérieur $dd^c u_{2,R} \wedge [\Gamma_a]$ est géométrique car le théorème de Fubini s'applique et l'intersection pluripotentialiste de deux courbes est géométrique. $\square$

Il est alors naturel de se poser la question suivante : sous quelles conditions le produit extérieur de deux courants uniformément laminaires est il admissible ? Il est aisé de construire des exemples de courants uniformément laminaires d'*auto-intersection* non admissible. C'est en effet le cas pour un courant d'intégration sur une courbe : dans des coordonnées $(z, w)$ adaptées la courbe a pour équation $\{z = 0\}$ et le courant associé a pour potentiel $\log|z|$. Pour obtenir des exemples diffus, il suffit de considérer des courants de la forme $\int [z = \alpha] d\mu(\alpha)$ où $\mu$ est une mesure d'énergie infinie (i.e. $\int u d\mu = -\infty$, $u$ étant un potentiel de $\mu$). En considérant une seconde mesure $\mu_2 \geq 0$ telle que $\int u d\mu_2 = -\infty$, on obtient un exemple de deux courants distincts dont on ne peut pas prendre le produit extérieur.

La proposition suivante dit que ces exemples sont essentiellement les seuls.

**Proposition 3.3.** *Soient $T_1$ et $T_2$ des courants uniformément laminaires diffus. Supposons que les feuilles des laminations sous jacentes $\mathcal{L}_1$ et $\mathcal{L}_2$ ne se coupent qu'en des points isolés. Alors le produit extérieur $T_1 \wedge T_2$ est localement admissible.*

Avant de commencer la preuve, rappelons une situation géométrique simple dans laquelle le produit extérieur de deux courants est admissible : c'est le cas si $\Omega \subset \mathbb{C}^2$ est un ouvert borné pseudoconvexe et $T_1$, $T_2$ sont des courants positifs fermés tels que $\operatorname{Supp} T_1 \cap \operatorname{Supp} T_2 \subset\subset \Omega$ (voir par exemple [De1, FS]).

**Preuve :** soit $x \in \operatorname{Supp} T_1 \cap \operatorname{Supp} T_2$ ; nous allons montrer qu'il y a un voisinage de $x$ dans lequel le produit extérieur $T_1 \wedge T_2$ est admissible. Soient respectivement $L_1(y)$ et $L_2(y)$ les feuilles des laminations associées à $T_1$ et $T_2$ passant par $y$. Soit $V \ni x$ tel que $L_k(x) \cap V$, $k = 1, 2$ soient des sous variétés de $V$ et $L_1(x) \cap L_2(x) \cap V = \{x\}$. Alors pour $V' \ni x$ suffisamment petit, on a

$$(4) \qquad \bigcup_{y \in V'} L_1(y) \cap L_2(y) \subset\subset V.$$

Comme dans la preuve du théorème précédent, on écrit $T_i = T_i' + T_i''$, où $T_i'$ est le courant subordonné à $T_i$ formé des feuilles issues de $V'$ ; en particulier $T_i'|_{V'} = T_i|_{V'}$. La propriété (4) assure que le produit $T_1' \wedge T_2'$ est admissible dans $V$. $\square$

## 4. Courants fortement approximables

On a vu à la section 2 que la représentation des courants laminaires pose problème en vue de l'interprétation géométrique de leur produit extérieur. Ces problèmes peuvent être surmontés si on dispose d'une suite de courbes holomorphes convergeant vers le courant de manière contrôlée. C'est le cas pour les courants adhérents aux courbes entières injectives [BLS, Ca], et pour les limites de certaines suites de diviseurs rationnels [Du1]. Rappelons que ceci entraîne la laminarité des courants invariants par les applications birationnelles de $\mathbb{P}^2$. Nous ne savons pas cependant si les courants laminaires construits par H. de Thélin [dT] –qui donne un critère local de laminarité– satisfont la définition 4.1 ci-dessous. Rappelons quelques éléments des constructions présentées dans ces travaux.

Soit $\Omega$ un ouvert de $\mathbb{C}^2$. On considère une suite de sous ensembles analytiques, éventuellement à bord, définis dans un voisinage de $\Omega$, d'aire $d_n$ tendant vers l'infini, et tels que $d_n^{-1}[C_n] \rightharpoonup T$. On veut obtenir les disques de $T$ comme graphes pour une projection linéaire $\pi$. Soit $L$ une base de la projection orthogonale à la direction de $\pi$ ; pour une subdivision $\mathcal{Q}$ de $L$ en carrés de taille $r$, on dira qu'une composante connexe de $\pi^{-1}(Q) \cap C_n$, $Q \in \mathcal{Q}$, est *bonne* si c'est un graphe au dessus de $Q$, *mauvaise* sinon. On retirera également en général un certain nombre de bonnes composantes afin d'assurer une propriété de compacité (critère de Montel ou borne sur l'aire). Dans les situations mentionnées plus haut, on a un contrôle en $O(d_n)$ du nombre total de mauvaises composantes, comptées avec multiplicité. Ainsi ces courants satisfont la définition 4.1 ci-dessous.

Dans la suite, nous considèrerons des subdivisions en carrés de taille $r$, images par une isométrie affine de la subdivision standard

$$\bigcup_{(j,k) \in \mathbb{Z}^2} \{z \in \mathbb{C}, \ jr < \Re(z) < (j+1)r, \ kr < \Im(z) < (k+1)r\};$$

on omettra en général la mention de la base $L$ de projection. On notera $\omega_1$ la restriction à $L$ de la forme Kählerienne standard de $\mathbb{C}^2$, et la notation $C$ désignera une "constante", pouvant varier de ligne en ligne, mais toujours indépendament de $n$ et $r$.

**Définition 4.1.** *Soit $T$ un courant positif fermé dans $\Omega$. $T$ est fortement approximable s'il existe une suite $(C_n)$ de sous ensembles analytiques définis dans un voisinage de $\Omega$, éventuellement à bord, telles que $d_n^{-1}[C_n] \rightharpoonup T$, au moins deux projections linéaires $\pi$, et une constante $C$, telles que si $\mathcal{Q}$ est une subdivision en carrés de taille $r$, et $T_{\mathcal{Q},n}$ désigne le courant formé des bonnes composantes de $C_n$ au dessus de $\mathcal{Q}$, normalisé par $d_n$, on ait l'estimation*

$$\text{(5)} \qquad \qquad \langle d_n^{-1}[C_n] - T_{\mathcal{Q},n}, \pi^* \omega_1 \rangle \leq Cr^2.$$

Un courant satisfaisant cette définition est laminaire (nous en redonnerons une preuve ci après à la proposition 4.4). Les courants laminaires construits dans l'article [Du1] vérifient cette définition dans tout ouvert de $\mathbb{P}^2$ ; nous les dirons "fortement approximables dans $\mathbb{P}^2$". Les courants adhérents aux courbes entières injectives [BLS, Ca], de même que les courants invariants des applications d'allure Hénon [Du2] sont également fortement approximables.

On a le théorème d'intersection géométrique (local) suivant.

**Théorème 4.2.** *Soient $T_1$ et $T_2$ des courants fortement approximables dans $\Omega$, de potentiel continu. Alors le produit extérieur $T_1 \wedge T_2$ est géométrique dans $\Omega$.*

On pourra noter que le théorème est également valable lorsqu'il existe un ouvert $\Omega' \subset \Omega$ dans lequel les courants $T_1$ et $T_2$ sont de potentiel continu et la masse de $T_1 \wedge T_2|_{\Omega'}$ soit totale dans $\Omega$. Ceci s'applique par exemple si $T_1$ et $T_2$ sont de potentiel borné, continu hors d'un fermé pluripolaire, ou encore si les potentiels de $T_1$ et $T_2$ sont continus hors d'un nombre fini de points ou l'un au plus des $T_i$ admet un nombre de Lelong –une situation analogue apparaît naturellement dans la dynamique de certaines applications birationnelles [Di2]).

**Corollaire 4.3.** *Soit $T$ un courant fortement approximable dans $\Omega$, et de potentiel continu. Alors $T \wedge T = 0$ dans $\Omega$. En particulier il existe des courants laminaires qui ne sont pas fortement approximables.*

En effet nous avons donné à la section 2 des exemples de courants laminaires de potentiel continu et d'auto-intersection strictement positive. On sait par ailleurs qu'un courant positif fermé de $\mathbb{P}^2$ tel que $T \wedge T$ soit admissible est d'auto-intersection strictement positive pour des raisons homologiques [FS]. On déduit en particulier de ce corollaire que les courants laminaires de potentiel continu dans $\mathbb{P}^2$ –c'est le cas en particulier pour l'exemple de Demailly 2.3– ne peuvent pas être fortement approximables.

Une autre conséquence du corollaire est le fait que pour le courant de Green $T^+$ d'un automorphisme polynomial de $\mathbb{C}^2$, qui est de potentiel continu dans $\mathbb{C}^2$, l'auto-intersection $T^+ \wedge T^+$ est concentrée au point d'indétermination $I^+$. Cet exemple montre également qu'une hypothèse sur le potentiel des courants est nécessaire dans le théorème 4.2. Le courant considéré ici admet un nombre de Lelong strictement positif au point $I^+$.

De manière plus générale, l'auto-intersection géométrique d'un courant laminaire est par définition toujours nulle, donc si $T$ est un courant laminaire de $\mathbb{P}^2$ tel que $T \wedge T$ soit admissible, ce produit, qui est de masse strictement positive, n'est jamais géométrique. Il serait intéressant d'en comprendre la structure.

Comme évoqué dans l'introduction, le théorème 4.2 donne une approche relativement directe pour certaines propriétés fines des automorphismes polynomiaux de $\mathbb{C}^2$. Nous renvoyons au panorama de Sibony [Si] pour un exposé de leurs principales propriétés. Soit par exemple $f$ un automorphisme de Hénon de $\mathbb{C}^2$, et $p$ un point périodique de type selle, que nous supposerons fixe sans perte de généralité. Alors les variétés stable et instable globales $W^s(p)$ et $W^u(p)$ sont des immersions injectives de $\mathbb{C}$ dans $\mathbb{C}^2$. On peut associer à de telles immersions des "courants d'intégration" fermés, qui sont des limites de courants d'intégration sur de grands disques, et sont laminaires d'après [BLS]. De plus leur structure laminaire est "limite" de celle des grands disques approximants. Si $T^s$ est associé à $W^s(p)$ et normalisé, alors Supp($T^s$) $\subset K^+$. D'après un théorème de Fornæss-Sibony, ceci implique $T^s = T^+$. On obtient de façon analogue $T^u = T^-$ pour un courant d'intégration sur la variété instable. Comme $T^+$ et $T^-$ sont de potentiel continu et $T^+ \wedge T^- > 0$, le théorème 4.2 s'applique et force les variétés stable et instable à se couper (transversalement). Ceci a des conséquences dynamiques (théorème homocline de Smale) : par exemple $p$ est accumulé par des points selles. Plus généralement ceci permet de retrouver directement tous les résultats du §9 de [BLS]. La même méthode s'applique au cas des applications birationnelles (très régulières) telles que les ensembles d'indétermination de $f$ et $f^{-1}$ soient dans l'ensemble de Fatou ("completely separating" selon la terminologie de Diller [Di1] §6).

Le théorème 4.2 peut être également utilisé comme résultat de convergence dans des situations où la théorie du pluripotentiel est délicate d'emploi. Nous illustrons ceci sur un exemple.

Soit $f$ un automorphisme polynomial de $\mathbb{C}^2$, $L$ et $L'$ deux droites complexes de $\mathbb{C}^2$. On a

$$\frac{1}{d^n}(f^n)_*[L] \rightharpoonup T^- \text{ et } \frac{1}{d^n}(f^n)^*[L] \rightharpoonup T^+.$$

On montre facilement à l'aide du théorème qu'on a la caractérisation suivante de la mesure $\mu = T^+ \wedge T^-$ :

$$\mu = \lim_{n \to \infty} [f^n(L) \cap f^n(L')].$$

Les méthodes usuelles de théorie du pluripotentiel ne permettent pas de démontrer un tel résultat.

Le reste de cette section sera consacré à la preuve du théorème 4.2. Nous commençons par montrer une version quantitative de la définition 2.1. Certains arguments peuvent être simplifiés dans le cas des courants fortement approximables dans $\mathbb{P}^2$, cependant nous préférons donner le résultat dans toute sa généralité. Le théorème étant local on pourra supposer que $\Omega$ est une boule.

**Proposition 4.4.** *Soit $T$ un courant laminaire fortement approximable dans $\Omega$. Fixons $\Omega' \subset\subset \Omega$. Alors pour des projections linéaires distinctes $\pi_1$, $\pi_2$ comme à la définition 4.1 et pour toutes subdivisions $\mathcal{S}_1$, $\mathcal{S}_2$ de $\mathbb{C}$ en carrés de taille $r$,*

$$\mathcal{Q} = \left\{ \pi_1^{-1}(s_1) \cap \pi_2^{-1}(s_2), (s_1, s_2) \in \mathcal{S}_1 \times \mathcal{S}_2 \right\}$$

*est une subdivision d'un voisinage de $\overline{\Omega'}$ en 4-cubes affines et il existe un courant $T_{\mathcal{Q}} \leq T$, uniformément laminaire dans chaque $Q \in \mathcal{Q}$, tel que*

(6)                                $$\mathbf{M}_{\Omega'}(T - T_{\mathcal{Q}}) \leq Cr^2,$$

*où $c$, $C$ sont indépendants de $r$ ($\mathbf{M}_{\Omega'}$ désigne la masse restreinte à $\Omega'$).*

**Preuve :** la preuve est similaire à celle de la proposition 3.4 de [Du1], la différence étant que l'on doit contrôler la masse simultanément dans deux directions. Un point délicat vient de ce qu'on n'a pas de contrôle sur les bonnes composantes à l'extérieur de $\Omega$.

Soient une suite $(C_n)$ ainsi que deux projections satisfaisant les hypothèses de la définition 4.1. Quitte à réduire légèrement $\Omega$, on supposera que les courants $d_n^{-1}[C_n]$ y sont de masse uniformément bornée. Considérons deux subdivisions $\mathcal{S}_i$ de taille $r$ et formons les courants $T_{\mathcal{S}_i,n}^i$ en retirant à $d_n^{-1}[C_n]$ les mauvaises composantes relatives à $\pi_i$ et $\mathcal{S}_i$, ainsi que les bonnes composantes dont la restriction à $\Omega$ est de grand volume. Le nombre de ces composantes est majoré par $O(d_n)$, car la masse de $d_n^{-1}[C_n]$ est $O(d_n)$ dans $\Omega$. Ainsi l'estimée (5) est satisfaite :

$$\left\langle d_n^{-1}[C_n] - T_{\mathcal{S}_i,n}^i, \pi_i^* \omega_1 \right\rangle \leq Cr^2, \ i = 1, 2.$$

Pour $i = 1, 2$, le courant $T_{\mathcal{S}_i,n}^i$ est le courant d'intégration sur une courbe à bord (non nécessairement connexe), dont le bord est inclus dans $\pi_i^{-1}(\partial \mathcal{S}_i)$,

$$T_{\mathcal{S}_i,n}^i = \frac{1}{d_n}[C_{\mathcal{S}_i,n}^i].$$

Considérons la subdivision $\mathcal{Q}$ comme définie dans l'énoncé de la proposition ; soit le courant

$$T_{\mathcal{Q},n} = \frac{1}{d_n}[C_{\mathcal{S}_1,n}^1 \cup C_{\mathcal{S}_2,n}^2] =: \frac{1}{d_n}[C_{\mathcal{Q},n}].$$

Il ne suffit de le considérer que pour les cubes $Q \in \mathcal{Q}$ tels que $Q \cap \Omega' \neq \emptyset$. Par construction $C_{\mathcal{Q},n}$ est une courbe holomorphe à bord, dont le bord est inclus dans $\partial \mathcal{Q}$ et comme pour $i = 1, 2$,

$$d_n^{-1}[C_n] - T_{\mathcal{Q},n} \leq d_n^{-1}[C_n] - T_{\mathcal{S}_i,n},$$

on en déduit

$$\langle d_n^{-1}[C_n] - T_{\mathcal{Q},n}, \ \pi_1^*\omega_1 + \pi_2^*\omega_1 \rangle \leq 2Cr^2;$$

il est par ailleurs à noter que $\pi_1^*\omega_1 + \pi_2^*\omega_1$ est une forme Kählerienne.

Il reste à voir qu'une sous-suite de $T_{\mathcal{Q},n}$ converge vers un courant $T_{\mathcal{Q}}$, uniformément laminaire dans chaque $Q \in \mathcal{Q}$. Soit $Q \in \mathcal{Q}$, décomposons $C_{\mathcal{Q},n} \cap Q$ en union disjointe

$$C_{\mathcal{Q},n} \cap Q = (C_{\mathcal{S}_1,n}^1 \cap Q) \cup \left[ (C_{\mathcal{S}_2,n}^2 \cap Q) \backslash (C_{\mathcal{S}_1,n}^1 \cap Q) \right],$$

ainsi que $T_{Q,n} := T_{\mathcal{Q},n}|_Q = T_{Q,n}^1 + T_{Q,n}^2$ selon cette décomposition.

Toutes les courbes de $T_{Q,n}^1$ sont des restrictions à $Q$ de sous variétés de $\pi_1^{-1}(s_1) \cap \Omega$, de volume majoré, et qui sont des graphes pour $\pi_1$. Ainsi une valeur d'adhérence d'une suite de composantes de $T_{Q,n}^1$ est un sous ensemble analytique de $Q$ qui est soit restriction d'un graphe pour $\pi_1$ soit inclus dans une fibre de $\pi_1$.

On en déduit que les valeurs d'adhérence de $T_{Q,n}^1$ sont des courants uniformément laminaires. En effet soit $V$ une limite d'une suite de composantes de $T_{Q,n_j}^1$. On déduit de la discussion précédente que $V$ est lisse. Fixons $x \in V$ et $\tau$ une petite transversale à $V$ en $x$. Toutes les composantes de $T_{Q,n_j}^1$ coupant $\tau$ pour $n_j$ assez grand sont des graphes au dessus d'un voisinage de $x$ dans $V$ car les composantes sont disjointes (lemme d'Hurwitz). De plus la masse transverse est majorée car $d_n^{-1}[C_n]$ est de masse bornée. On conclut en utilisant les arguments usuels : voir par exemple [Du1] proposition 3.4 ou [BS5].

On fait le même raisonnement pour $T_{Q,n}^2$, et on considère une valeur d'adhérence $T_Q = T_Q^1 + T_Q^2$. Si $x \in \text{Supp}(T_Q^1) \cap \text{Supp}(T_Q^2)$, et $\Delta^1$, $\Delta^2$, sont les disques correspondants passant par $x$, alors $\Delta^1 = \Delta^2$. En effet si $x$ était un point isolé de $\Delta^1 \cap \Delta^2$, les disques approximants respectifs $\Delta_{n_j}^1$ et $\Delta_{n_j}^2$ des $T_{Q,n_j}^i$ devraient se couper, ce qui est impossible puisque ces disques sont des bonnes composantes de $C_n$. Ceci montre que $T_Q$ est uniformément laminaire. □

**Preuve du théorème 4.2 :** soit $\nu$ la mesure $T_1 \wedge T_2$, nous supposerons que $\|\nu\| \leq 1$ dans $\Omega$. Nous allons construire une suite croissante de mesures d'intersection géométrique $\nu_{\mathcal{Q}} \leq \nu$ telles que $\mathbf{M}(\nu - \nu_{\mathcal{Q}}) \to 0$ quand le pas de la subdivision tend vers 0.

On dispose pour $k = 1, 2$ de subdivisions –dans le cas de courants fortement approximables dans $\mathbb{P}^2$ on peut choisir une subdivision commune aux deux courants– en 4-cubes affines $\mathcal{Q}^k$ de taille comparable à $r$ (vues comme restrictions à $\Omega$ de subdivisions de $\mathbb{C}^2$) telles que les conclusions de la proposition précédente soient satisfaites pour $T_k$. Étant donné que les subdivisions $\mathcal{S}$ en carrés de taille $r$ ne sont pas fixées, pour $Z_k \in \mathbb{C}^2$, ceci vaut également pour les subdivisions $Z_k + \mathcal{Q}^k$. Sans perte de généralité on écrira dorénavant $\Omega$ pour $\Omega'$.

Posons $\mathcal{Q} = \mathcal{Q}^k$ pour $k$ valant 1 ou 2, et $\mathcal{Q}_Z$ la subdivision translatée $\mathcal{Q} + Z$. Si $Z$ est dans le réseau $L(\mathcal{Q})$ des sommets de $\mathcal{Q}$, alors $\mathcal{Q}_Z = \mathcal{Q}$. Le nombre de cubes rencontrant $\Omega$ est $\text{Card}_\Omega(\mathcal{Q}_Z) \leq C_0/r^4$, où $C_0$ est indépendante de $Z$.

Soit un nombre réel $\lambda < 1$, et pour $Q \in \mathcal{Q}$, soit $Q^\lambda$ l'homothétie de $Q$ de rapport $\lambda$ par rapport à son centre. Le volume de Lebesgue $\text{Leb}(Q^\lambda)$ est $\lambda^4 \text{Leb}(Q)$. Soit $\mathcal{Q}^\lambda := \cup_{Q \in \mathcal{Q}} Q^\lambda$,

$(\mathcal{Q}^\lambda)^c := \mathbb{C}^2 \backslash \mathcal{Q}^\lambda$,

$$\mathrm{Leb}(\Omega \cap (\mathcal{Q}^\lambda)^c) \le \mathrm{Card}_\Omega(\mathcal{Q})(1 - \lambda^4) \mathrm{Leb}(Q) \le C(1 - \lambda^4).$$

**Lemme 4.5.** *Il existe $Z \in \mathbb{C}^2$ tel que $\nu(\Omega \cap (\mathcal{Q}_Z^\lambda)^c) < 2(1 - \lambda^4)$.*

**Preuve :** le lemme est une conséquence simple de l'invariance de la mesure de Lebesgue par translation et du théorème de Fubini. Par périodicité il suffit de considérer $Z \in \overline{Q}$ pour un certain $Q \in \mathcal{Q}$. On a

$$\int_{\overline{Q}} \nu(\Omega \cap (\mathcal{Q}_Z^\lambda)^c) \frac{d\,\mathrm{Leb}(Z)}{\mathrm{Leb}(Q)} = \int_{\overline{Q}} \int_\Omega \mathbf{1}_{Z + (\mathcal{Q}^\lambda)^c}(y) d\nu(y) \frac{d\,\mathrm{Leb}(Z)}{\mathrm{Leb}(Q)}$$

et $y \in Z + (\mathcal{Q}^\lambda)^c$ si et seulement si $Z \in y - (\mathcal{Q}^\lambda)^c$. Donc $\mathbf{1}_{Z + (\mathcal{Q}^\lambda)^c}(y) = \mathbf{1}_{y - (\mathcal{Q}^\lambda)^c}(Z)$. L'intégrale précédente vaut donc

$$\int_\Omega \left( \int_{\overline{Q}} \mathbf{1}_{y - (\mathcal{Q}^\lambda)^c}(Z) \frac{d\,\mathrm{Leb}(Z)}{\mathrm{Leb}(Q)} \right) d\nu(y) = \int_\Omega \frac{1}{\mathrm{Leb}(Q)} \mathrm{Leb}((y - (\mathcal{Q}^\lambda)^c) \cap Q) d\nu(y),$$

et pour tout $y \in \mathbb{C}^2$, $\mathrm{Leb}((y - (\mathcal{Q}^\lambda)^c) \cap Q) = \mathrm{Leb}(Q \backslash \mathcal{Q}^\lambda) = (1 - \lambda^4) \mathrm{Leb}(Q)$, par invariance par translation de la mesure de Lebesgue sur le tore $\mathbb{C}^2/L(\mathcal{Q})$. D'où

$$\int_{\overline{Q}} \nu(\Omega \cap (\mathcal{Q}_Z^\lambda)^c) \frac{d\,\mathrm{Leb}(Z)}{\mathrm{Leb}(Q)} \le (1 - \lambda^4)$$

et le lemme est démontré.                                                              □

On fixe donc $Z_1$ et $Z_2$ tels que la conclusion du lemme soit respectivement satisfaite pour $\mathcal{Q}^1$ et $\mathcal{Q}^2$, et par abus de notation, on renomme $\mathcal{Q}^k$ la subdivision translatée $Z_k + \mathcal{Q}^k$. On applique la proposition 4.4 aux courants $T_k$, et ainsi on obtient des courants $T_{k,\mathcal{Q}^k} \le T_k$, uniformément laminaires dans chaque $Q \in \mathcal{Q}^k$, et tels que $\mathbf{M}(T_k - T_{k,\mathcal{Q}^k}) \le Cr^2$.

Considérons la subdivision $\mathcal{Q} = \mathcal{Q}^1 \wedge \mathcal{Q}^2$ formée des $Q^1 \cap Q^2$, $(Q^1, Q^2) \in \mathcal{Q}^1 \times \mathcal{Q}^2$. Soient $Q \in \mathcal{Q}$ et $T_{k,Q} = T_{k,\mathcal{Q}^k}|_Q$ ; $T_{k,Q} \le T_k$ donc le produit $T_{1,Q} \wedge T_{2,Q}$ est admissible dans $Q$ (les $T_{k,Q}$ sont en fait à potentiel continu d'après [BLS] lemme 8.2), et est géométrique par le théorème 3.1. Soit

$$\nu_\mathcal{Q} = T_{1,\mathcal{Q}} \wedge T_{2,\mathcal{Q}} = \sum_{Q \in \mathcal{Q}} T_{1,Q} \wedge T_{2,Q} = \sum_{Q \in \mathcal{Q}} T_{1,Q} \dot{\wedge} T_{2,Q};$$

on a $\nu_\mathcal{Q} \le \nu$ et il reste à estimer $\mathbf{M}(\nu - \nu_\mathcal{Q})$.

On pose pour $\lambda < 1$,

$$\mathcal{Q}^\lambda = \bigcup_{(Q^1, Q^2) \in \mathcal{Q}^1 \times \mathcal{Q}^2} (Q^1)^\lambda \cap (Q^2)^\lambda,$$

(on peut avoir $Q_1 \cap Q_2 \ne \emptyset$ et $(Q^1)^\lambda \cap (Q^2)^\lambda = \emptyset$). La masse se décompose alors en

$$\mathbf{M}_\Omega(\nu - \nu_\mathcal{Q}) = \mathbf{M}_{\Omega \cap \mathcal{Q}^\lambda}(\nu - \nu_\mathcal{Q}) + \mathbf{M}_{\Omega \cap (\mathcal{Q}^\lambda)^c}(\nu - \nu_\mathcal{Q}).$$

Le deuxième terme est majoré par $4(1 - \lambda^4)$ d'après le lemme précédent. On fixe $\lambda$ tel que $2(1 - \lambda^4) < \varepsilon/2$ et il reste à majorer $\mathbf{M}_{\Omega \cap \mathcal{Q}^\lambda}(\nu - \nu_\mathcal{Q})$.

Soit $Q \in \mathcal{Q}$, alors

$$T_1 \wedge T_2 - T_{1,Q} \wedge T_{2,Q} = (T_1 - T_{1,Q}) \wedge T_2 + T_{1,Q} \wedge (T_2 - T_{2,Q})$$
$$\le (T_1 - T_{1,Q}) \wedge T_2 + T_1 \wedge (T_2 - T_{2,Q}).$$

Par symétrie il suffit de considérer le premier terme. Soit $u_2$ un potentiel de $T_2$ dans $\Omega$ ; on va construire une fonction $Q$ dont la restriction à chaque $Q$ est une fonction plateau, $\chi|_Q \in C_0^\infty(Q)$, $\chi \geq 0$, $\chi = 1$ au voisinage de $Q^\lambda$, et $\|dd^c\chi\|_{L^\infty} = O(1/r^2)$.

Dans un carré de taille $r$ dans $\mathbb{C}$ on construit aisément une telle fonction, avec $\|d\chi_0\|_{L^\infty} \leq C/(1-\lambda)^2 r^2$. D'où par translation une fonction $\chi_0$ ayant les propriétés souhaitées dans une subdivision en carrés de taille $r$ dans $\mathbb{C}$. Si maintenant $\pi_i$, $i = 1, \ldots, 4$, sont les quatre projections linéaires correspondant aux deux subdivisions $\mathcal{Q}^1$ et $\mathcal{Q}^2$, il suffit de poser

$$\chi = \chi_0 \circ \pi_1 \cdots \chi_0 \circ \pi_4.$$

Le paramètre $\lambda$ étant fixé, $\|dd^c\chi\|_{L^\infty} = O(\frac{1}{r^2})$. Par ailleurs

$$\mathbf{M}_{Q^\lambda}((T_1 - T_{1,Q}) \wedge T_2) \leq \int \chi(T_1 - T_{1,Q}) \wedge T_2) = \int u_2 dd^c\chi \wedge (T_1 - T_{1,Q}),$$

et si $c_Q$ est un point dans $Q$, on remarque que

$$\int u_2 dd^c\chi \wedge (T_1 - T_{1,Q}) = \int (u_2 - u_2(c_Q)) dd^c\chi \wedge (T_1 - T_{1,Q}).$$

Donc

$$\int u_2 dd^c\chi \wedge (T_1 - T_{1,Q}) \leq C\omega(u_2, r)\mathbf{M}(T_1 - T_{1,Q}) \|dd^c\chi\|_{L^\infty} \leq C\omega(u_2, r)\frac{1}{r^2}\mathbf{M}(T_1 - T_{1,Q}),$$

où $\omega(u_2, r)$ est le module de continuité de $u_2$ de rayon $r$. En prenant la somme sur $Q \in \mathcal{Q}$, et en utilisant (6) on obtient

$$\sum_{Q \in \mathcal{Q}} \mathbf{M}_{Q^\lambda}((T_1 - T_{1,Q}) \wedge T_2) \leq C\omega(u_2, r),$$

soit

$$\mathbf{M}_{\Omega \cap \mathcal{Q}^\lambda}(\nu - \nu_{\mathcal{Q}}) \leq C(\omega(u_1, r) + \omega(u_2, r)).$$

Donc si $r$ est suffisamment petit, $\mathbf{M}_\Omega(\nu - \nu_{\mathcal{Q}}) < \varepsilon$.

En reproduisant inductivement le même procédé dans chaque $Q \in \mathcal{Q}$, on obtient bien une suite croissante de mesures $\nu_{\mathcal{Q}}$                          $\square$

**Remarque 4.6.** On déduit de la preuve que si l'un de $T_1$ ou $T_2$, par exemple $T_1$, est uniformément laminaire, alors les termes correspondant à $(T_1 - T_{1,Q})$ sont nuls. En particulier dans ce cas l'hypothèse de potentiel continu pour $T_2$ est inutile. Nous avons donc montré le résutat suivant : si $T_1$ est un courant uniformément laminaire et à potentiel continu dans $\Omega$, et $T_2$ est un courant fortement approximable quelconque, alors l'intersection $T_1 \wedge T_2$ est géométrique dans $\Omega$.

**Remarque 4.7.** La preuve donne une estimation explicite de la masse de $\nu - \nu_{\mathcal{Q}}$. Pour le choix $1 - \lambda = (\omega(u_1, r) + \omega(u_2, r))^{1/3}$ on obtient une majoration de $\mathbf{M}(\nu - \nu_{\mathcal{Q}})$ en $(\omega(u_1, r) + \omega(u_2, r))^{1/3}$.

Mathématique, Bâtiment 425, Université de Paris Sud, 91405 Orsay cedex, France.
`Romain.Dujardin@math.u-psud.fr`